\documentclass[11pt,reqno]{amsart}
\usepackage[utf8]{inputenc}
\usepackage[english]{babel}
\usepackage{amsmath,amsfonts,amssymb,amsthm}
\usepackage{hyperref}

\makeatletter
\@namedef{subjclassname@2020}{\textup{2020} Mathematics Subject Classification}
\makeatother

\newcommand{\Z}{{\mathbb Z}}

\newcommand{\length}{{\rm length}}

\newcommand{\ent}{{\rm ent}}
\newcommand{\dist}{{\rm dist}}
\newcommand{\vol}{\operatorname{vol}}
\newcommand{\area}{\operatorname{area}}

\numberwithin{equation}{section}

\newtheorem{theorem}{Theorem}[section]
\newtheorem{proposition}[theorem]{Proposition}
\newtheorem{corollary}[theorem]{Corollary}
\newtheorem{lemma}[theorem]{Lemma}

\theoremstyle{definition}
\newtheorem{definition}[theorem]{Definition}
\newtheorem{example}[theorem]{Example}
\newtheorem{remark}[theorem]{Remark}

\begin{document}

\title{Macroscopic scalar curvature bounds on graphs}

\author[S.~Karam]{Steve Karam}
\address{College of Engineering and Technology, American University of the Middle East, Egaila 54200, Kuwait}
\email{steve.karam@aum.edu.kw}

\author[S.~Sabourau]{St\'ephane Sabourau}
\address{\parbox{\linewidth}{Univ Paris Est Creteil, CNRS, LAMA, F-94010 Creteil, France \\
    Univ Gustave Eiffel, LAMA, F-77447 Marne-la-Vall\'ee, France}}
\email{stephane.sabourau@u-pec.fr}


\begin{abstract}
  We establish sharp lower bounds for the lengths of balls in the universal covers of finite metric graphs.
  For a finite connected $d$-regular graph, with $d \geq 3$ and average edge length at most one, we show that, at every radius, the average length of balls centered at a vertex of the universal cover is at least that of a ball of the same radius in the equilateral $d$-regular tree.
  We extend this comparison to arbitrary finite connected graphs whose total length does not exceed that of an equilateral trivalent graph with the same cyclomatic number.
  In particular, we show that the maximal length of a ball centered at a vertex of the universal cover is bounded below by that of a ball of the same radius in the equilateral trivalent tree.
  This answers a question of the first author.
  Furthermore, equality at an integer radius forces all edges to have unit length in the regular case and, in the general case, forces the graph to be an equilateral trivalent graph after suppressing vertices of degree two.
  In contrast, rigidity fails at every noninteger radius.
\end{abstract}

\maketitle

\section{Introduction}

Given a closed Riemannian manifold~$(M,g)$, denote by~$(\tilde{M},\tilde{g})$ its universal cover endowed with the lifted metric.
Following a viewpoint going back to Gromov~\cite{gro86}, the volumes of balls of a fixed radius~$r$ in~$\tilde M$ can be regarded as a macroscopic analogue of scalar curvature at scale~$r$; see also~\cite{Gut10}.
This viewpoint is motivated by the fact that the scalar curvature appears in the asymptotic expansion of the volumes of small geodesic balls.

In this context, Guth~\cite[Conjecture~2]{Gut11} formulated the following sharp conjecture.
Let $(M,g_0)$ be a closed hyperbolic $n$-manifold, and let $g$ be another Riemannian metric on~$M$ with $\vol_g(M) < \vol_{g_0}(M)$.
Then, for every $r>0$,
\begin{equation} \label{eq:Guth}
  \max_{\tilde x \in \tilde M} \vol_{\tilde g} B(\tilde x,r) > \vol B_{\mathbb{H}^n}(r),
\end{equation}
where $B_{\mathbb{H}^n}(r)$ denotes a ball of radius~$r$ in the hyperbolic space~$\mathbb{H}^n$.
For sufficiently large~$r$, depending on~$g$, this inequality follows from the minimal volume entropy theorem of Besson, Courtois and Gallot~\cite{BCG}, and from Katok's theorem~\cite{katok} in dimension two.

In~\cite{Gut11}, Guth proved~\eqref{eq:Guth} for~$r=1$ under a sufficiently small volume ratio, namely $\vol_g(M) < \delta_n \vol_{g_0}(M)$ for some $\delta_n>0$ depending only on~$n$; see~\cite{BS21} and~\cite{A25} for further extensions.
He also asked whether this volume condition implies the inequality~\eqref{eq:Guth} simultaneously for all $r \geq 1$.
The first author~\cite{karam} answered this question affirmatively for surfaces.
In arbitrary dimensions, volume comparisons for all $r \geq 1$ were subsequently obtained under additional hypotheses: by Balacheff and the first author~\cite{BK19}, assuming that $g$ has negative sectional curvature and a sufficiently small volume ratio; by the second author~\cite{Sab22} for metrics of sufficiently small volume; and by Alpert~\cite{Alp26}, under the condition $\vol_g(M)^2<\delta'_n \vol_{g_0}(M)$ for some $\delta'_n>0$.
More recently, Zhang~\cite{zhang} obtained such a comparison under a sufficiently small volume ratio, thereby answering Guth's question in every dimension.

Returning to Guth's sharp conjecture, Chodosh and Vikman~\cite{CV} recently proved it for surfaces.
More precisely, let $\Sigma$ be a closed oriented surface of genus at least two, endowed with a Riemannian metric~$g$ of the same area as a hyperbolic metric on~$\Sigma$.
Then, for every $r>0$,
\begin{equation} \label{eq:CV}
  \max_{\tilde x \in \tilde \Sigma} \area B(\tilde x,r) \geq \area B_{\mathbb{H}^2}(r)
\end{equation}
with equality for a given~$r$ if and only if $g$ is hyperbolic.
Their proof also establishes a sharp lower bound for a suitable average of the areas of these balls.

\medskip

In this note, we consider the analogous problem for finite metric graphs, with regular trees playing the role of the hyperbolic plane.
The asymptotic counterpart of this problem is the minimization of volume entropy, which measures the exponential
growth rate of the lengths of balls in the universal cover; see Corollary~\ref{coro:entropy} below.
The minimal volume entropy of finite graphs was determined independently by Lim~\cite{Lim08} and by Kapovich and Nagnibeda~\cite{KN07}; it can also be deduced from the earlier work of Rivin~\cite{Rivin10}.
In particular, among finite connected graphs with a fixed cyclomatic number $c \geq 2$ and total length $3(c-1)$, the volume entropy is minimal precisely for trivalent graphs with unit edge lengths.
See also~\cite{robert} for graphs with highly transitive automorphism groups, \cite{McMullen15} for another approach, \cite{sab} for an elementary proof, and~\cite{bal,KL20} for related results.

Motivated by Guth's question, the first author~\cite{karam} asked whether there exists a universal constant $\kappa>0$ such that every finite connected metric graph whose total length is less than $\kappa$ times that of an equilateral trivalent graph with the same cyclomatic number has the following property: for every radius $r>0$, the maximal length of a ball of radius~$r$ in its universal cover is greater than or equal to the maximal length of a ball of the same radius in the equilateral trivalent tree.
We answer this question with the optimal constant~$\kappa=1$ in Theorem~\ref{theo:nonreg} below.

\medskip

We first establish a sharp comparison at every radius for regular graphs, with regular trees playing the role of the
hyperbolic plane.
This is a graph analogue of~\cite[Theorem~4]{CV}.
As in~\cite{CV}, we obtain a lower bound for an average of ball lengths, where the average is taken over the vertices of the graph.

\medskip

We begin by recalling some basic definitions about graphs.

\begin{definition}
  Let $X$ be a finite connected graph, where loops and multiple edges are allowed.
  Denote by $V(X)$ and~$E(X)$ the sets of vertices and of unoriented edges of~$X$.
  By a slight abuse of notation, $V(X)$ and~$E(X)$ also denote their numbers of elements.

  A graph~$X$ is \emph{$d$-regular} if every vertex of~$X$ has degree~$d$ (with loops counted twice).
  A \emph{length metric} on~$X$ assigns a positive length~$\ell(e)$ to each edge~$e$.
  Denote by~$\length(X)$ the sum of these lengths and by~$\tilde X$ the universal cover of~$X$ with the lifted metric.

  Let $T_d$ be the $d$-regular tree with all edges of unit length.
  For every $r>0$, denote by
  \begin{equation} \label{eq:Td}
    V_{T_d}(r) = \length \, B_{T_d}(x,r)
  \end{equation}
  the length of an $r$-ball centered at a vertex~$x$ of~$T_d$.
  For $d \geq 3$, we have
  \[
    V_{T_d}(r) = \frac{d(d-1)^{\lfloor r \rfloor}-d}{d-2} + d(d-1)^{\lfloor r \rfloor} (r- \lfloor r \rfloor).
  \]
  This follows by summing the $d(d-1)^{k-1}$ edges between the vertices at distances~$k-1$ and~$k$ from~$x$.
\end{definition}

Our first result is the following sharp comparison for regular graphs.

\begin{theorem} \label{theo:main}
  Let $X$ be a finite connected $d$-regular graph with~$d \geq 3$, endowed with a length metric satisfying $\length(X) \leq E(X)$.
  Then, for every $r>0$,
  \begin{equation} \label{eq:main-average}
    \frac{1}{V(X)} \sum_{x \in V(X)} \length \, B(\tilde{x},r) \geq V_{T_d}(r)
  \end{equation}
  where $\tilde{x}$ is any lift of~$x \in X$ to the universal cover~$\tilde{X}$.

  In particular,
  \begin{equation} \label{eq:main}
    \max_{x \in V(X)} \length \, B(\tilde{x},r) \geq V_{T_d}(r).
  \end{equation}

  Moreover, equality holds in either~\eqref{eq:main-average} or~\eqref{eq:main} at some integer~$r=n \geq 1$ if and only if every edge of~$X$ has unit length.
\end{theorem}

Both inequalities are optimal.
Indeed, equality holds for every radius when all edges of~$X$ have unit length, since its universal cover is then isometric to~$T_d$.

\medskip

The rigidity statement at integer radii contrasts with the behavior at the other radii.
In Proposition~\ref{prop:equality}, we characterize all radii at which equality holds in~\eqref{eq:main-average}, in terms of the shortest and longest edge-paths with a prescribed number of edges.
As a consequence, we show that at every noninteger radius, equality in~\eqref{eq:main-average} or in~\eqref{eq:main} does not force the edges to have unit length; see Example~\ref{ex:noninteger}.
We obtain a stronger non-rigidity phenomenon in Proposition~\ref{prop:theta}: there exist non-equilateral theta graphs for which equality holds in~\eqref{eq:main} on a nontrivial interval of radii, even when the maximum is taken over all centers in the universal cover.
Thus, unlike the comparison~\eqref{eq:CV} for surfaces, equality at a single radius does not generally imply rigidity.

\medskip

We now extend the comparison to arbitrary finite connected graphs, without any regularity assumption.
The resulting bound depends only on the cyclomatic number and is sharp.

\medskip

For this purpose, we recall the following definitions.

\begin{definition}
  Let $X$ be a finite connected graph, endowed with a length metric.
  The \emph{cyclomatic number} of~$X$ is defined as
  \[
    c(X)=E(X)-V(X)+1.
  \]
  It is the rank of~$H_1(X;\Z)$, or equivalently, the number of independent cycles of $X$.

  Suppose $c(X) \geq 2$.
  By successively removing vertices of degree one and their incident edges, and then suppressing vertices of degree two while adding the lengths of their incident edges, we obtain a metric graph~$K$, called the \emph{core} of~$X$.
  It has the same cyclomatic number as~$X$, and all its vertices have degree at least three.
  In particular,
  \[
    \sum_{x \in V(K)} \left( \deg_K(x)-2 \right) = 2E(K)-2V(K) = 2(c(X)-1).
  \]
\end{definition}

Our second result is the following sharp comparison for arbitrary graphs.

\begin{theorem} \label{theo:nonreg}
  Let $X$ be a finite connected graph with $c=c(X) \geq 2$, endowed with a length metric satisfying
  \[
    \length(X) \leq 3(c-1).
  \]
  Let $K$ be the core of~$X$.
  Then, for every $r>0$,
  \begin{equation} \label{eq:nonreg-average}
    \frac{1}{2(c-1)} \sum_{x \in V(K)} \left( \deg_K(x)-2 \right)  \length \, B_{\tilde{K}}(\tilde{x},r) \geq V_{T_3}(r)
  \end{equation}
  where $\tilde{x}$ is any lift of~$x$ to~$\tilde{K}$.

  In particular,
  \begin{equation} \label{eq:nonreg}
    \max_{x \in V(K)} \length \, B_{\tilde{X}}(\tilde{x},r) \geq V_{T_3}(r)
  \end{equation}
  where every vertex~$x$ of~$K$ is regarded as a vertex of~$X$, and $\tilde x$ is any lift of~$x$ to~$\tilde X$.

  Moreover, equality holds in either~\eqref{eq:nonreg-average} or~\eqref{eq:nonreg} at some integer~$r=n \geq 1$ if and only if $X$ is a trivalent graph with all edges of unit length, after suppressing its vertices of degree two.
\end{theorem}

Both inequalities are sharp, with equality for trivalent graphs~$X$ whose edges all have unit length.
In this case, $\length(X)=E(X)=3(c-1)$, the core of~$X$ is $X$ itself, its universal cover~$\tilde{X}$ is~$T_3$, and $\length \, B_{\tilde{X}}(\tilde{x},r) = V_{T_3}(r)$.
In particular, since balls in~$T_3$ are longest when centered at vertices, the theorem gives the optimal constant in the ball growth question considered in the introduction of~\cite{karam}.

\medskip

We emphasize that although rigidity fails at noninteger radii even for trivalent graphs, equality at an integer radius forces not only the metric to be equilateral but also the underlying graph to be trivalent after suppressing its vertices of degree two.

\medskip

As a consequence of our comparison theorems, we recover classical lower bounds for the minimal volume entropy of graphs.
Recall that the \emph{volume entropy} of a finite connected metric graph~$X$ is defined as the exponential growth rate of the length of balls in the universal cover.
That is,
\[
  \ent(X) = \lim_{r \to \infty} \frac{1}{r} \log \left( \length \, B_{\tilde{X}}(\tilde x,\,r) \right)
\]
where the limit exists and is independent of the center~$\tilde x$; see~\cite{man}.

Since $V_{T_d}(r)$ grows exponentially at rate~$\log(d-1)$, Theorem~\ref{theo:main} and Theorem~\ref{theo:nonreg} immediately yield the following classical result.

\begin{corollary} \label{coro:entropy}
  Let $X$ be a finite connected graph with cyclomatic number~$c=c(X) \geq 2$, endowed with a length metric.
  \begin{enumerate}
    \item If $X$ is $d$-regular with $d \geq 3$, and $\length(X) \leq E(X)$, then
          \[
            \ent(X) \geq \log(d-1).
          \]
    \item If $\length(X) \leq 3(c-1)$, then
          \[
            \ent(X) \geq \log(2).
          \]
  \end{enumerate}
  Equality is attained by regular graphs with all edges of unit length, of degree~$d$ in the first case and degree three in the second.
\end{corollary}

These entropy inequalities, including the uniqueness of their minimizers, are known; see~\cite{Rivin10,Lim08,KN07,McMullen15,sab} and the references therein.
Our results strengthen these asymptotic inequalities by providing sharp comparisons at every radius, although our method does
not recover the uniqueness of the entropy minimizers.

\medskip

The proof of Theorem~\ref{theo:main} follows the elementary approach developed in~\cite{sab}.
It relies on an expression of the length of balls in the universal cover in terms of edge-paths, an averaging argument for these edge-paths using the fact that every edge has unit length on average, and Jensen's convexity inequality.
Theorem~\ref{theo:nonreg} is then obtained by reducing the general case to the trivalent case, replacing vertices of higher degree
by trees whose internal edges have length zero.

\medskip

The article is organized as follows.
Section~\ref{sec:prelim} introduces the notation and establishes basic formulas for the lengths of balls in universal covers of metric graphs.
In Section~\ref{sec:proof13}, we prove Theorem~\ref{theo:main}, characterize its equality case, and show that rigidity fails at every noninteger radius.
In Section~\ref{sec:nonreg}, we deduce Theorem~\ref{theo:nonreg}.
Finally, Section~\ref{sec:theta} studies theta graphs and establishes a stronger non-rigidity result in which the center of the balls may vary over the entire universal cover.

\medskip

\emph{AI disclosure.} We used AI tools to revise and improve the presentation of the article.

\section{Formulas for ball length} \label{sec:prelim}

In this section, we introduce the notation for oriented edge-paths and establish formulas for the length of balls in universal covers of metric graphs.
We allow edges of length zero, provided that the subgraph formed by these edges contains no cycle.
In this case, the universal cover is equipped with the pseudometric induced by the edge lengths, and the length of a ball is defined as the sum of the lengths of its intersections with the edges.

\begin{definition}\label{def:path}
  Let $X$ be a finite connected graph, possibly with loops and multiple edges.
  Each edge of~$X$ gives two \emph{oriented edges}, one for each direction of travel.
  For an oriented edge~$a$, we denote by~$a^-$ and~$a^+$ its initial and terminal vertices, and by~$\bar a$ the same edge with the opposite orientation.
  For a loop, $a^-=a^+$.
  An \emph{oriented edge-path} of~$X$ of combinatorial length~$k \geq 1$ is a sequence $\gamma = a_1 \cdots a_k$ of oriented edges such that
  \[
    a_i^+=a_{i+1}^- \qquad \textrm{ and } \qquad a_{i+1} \neq \bar{a}_i
  \]
  for every $i=1,\dots,k-1$.
  In particular, edge-paths do not backtrack.

  When $X$ is endowed with edge lengths, we denote by~$\ell(a)$ the length of the underlying edge.
  The length of~$\gamma$ is then $\ell(\gamma)=\ell(a_1)+\dots+\ell(a_k)$.

  For every vertex $x \in X$, denote by $\emptyset_x$ the empty edge-path at $x$, of length zero.
  If $\gamma$ is an oriented edge-path ending at the starting point~$a^-$ of an oriented edge~$a$, we write $\gamma a$ for their concatenation, provided that no backtracking occurs.
  This convention also applies when $\gamma=\emptyset_x$.

  For $k \geq 0$, denote by $\mathcal{P}_k$ the set of oriented edge-paths of~$X$ of combinatorial length~$k$, and by $\mathcal{P}_k(x)$ the set of those starting at a vertex~$x$.
\end{definition}

The following observation will be used in the proofs of Lemmas~\ref{lem:ball_paths} and~\ref{lem:ball_all}.

\begin{remark} \label{rem:correspondence}
  Let $x$ be a vertex of~$X$.
  Fix a lift~$\tilde{x}$ in the universal cover~$\tilde{X}$ of~$X$.
  There is a one-to-one correspondence between the vertices~$\tilde{y} \in \tilde{X}$ and the edge-paths $\gamma \in \bigsqcup_{k \geq 0} \mathcal{P}_k(x)$.
  More precisely, the unique oriented edge-path from~$\tilde{x}$ to~$\tilde{y}$ projects to an oriented edge-path of~$X$ starting at~$x$.
  Conversely, every oriented edge-path of~$X$ starting at~$x$ lifts to a unique oriented edge-path starting at~$\tilde{x}$.
\end{remark}

\begin{example}
  Let $X$ be a finite connected $d$-regular graph, with $d\geq3$.
  For every vertex $x$ of~$X$ and every integer $k \geq 1$, there are exactly
  \begin{equation} \label{eq:Pkx}
    | \mathcal{P}_k(x) | = d(d-1)^{k-1}
  \end{equation}
  oriented edge-paths of combinatorial length~$k$ starting at~$x$.
  Therefore,
  \begin{equation} \label{eq:Pk}
    | \mathcal{P}_k |=d(d-1)^{k-1} V(X).
  \end{equation}
\end{example}

The following elementary averaging lemma appears in~\mbox{\cite[Lemma~7]{sab}}.
We reproduce its proof for completeness.

\begin{lemma} \label{lem:w-comb}
  Let $X$ be a finite connected $d$-regular graph endowed with nonnegative edge lengths.
  Then, for every integer~$k \geq 1$,
  \[
    \frac{1}{|\mathcal{P}_k|} \sum_{\gamma \in \mathcal{P}_k} \ell(\gamma) = k \, \frac{\length(X)}{E(X)}.
  \]
\end{lemma}

\begin{proof}
  Let $e_i$ be the unoriented edges of~$X$, with $1 \leq i \leq E(X)$.
  Every unoriented edge occurs exactly $2k \, (d-1)^{k-1}$ times among all oriented edge-paths of combinatorial length~$k$, counting multiplicities (the factor~$2$ accounts for the orientation of the edge-paths and the factor~$k$ for the position of the edge along the path).
  Therefore, the average length~$\ell(\gamma)$ of oriented edge-paths $\gamma \in \mathcal{P}_k$ of combinatorial length~$k$ is equal to
  \[
    \frac{1}{|\mathcal{P}_k|}  \sum_{\gamma \in \mathcal{P}_k} \ell(\gamma) = \frac{2 k \, (d-1)^{k-1}}{d(d-1)^{k-1} \, V(X)} \sum_{i=1}^{E(X)} \ell(e_i) = \frac{2k \, \length(X)}{d \, V(X)} = k \, \frac{\length(X)}{E(X)}.
  \]
  To derive these equalities, recall that $|\mathcal{P}_k| = d(d-1)^{k-1} \, V(X)$, $\sum_{i=1}^{E(X)} \ell(e_i) = \length(X)$ and $2 \, E(X)=d \, V(X)$.
\end{proof}

We next express the length of a ball in the universal cover in terms of oriented edge-paths.
Here and below, we write $s_+=\max\{s,0\}$.

\begin{lemma} \label{lem:ball_paths}
  Let $X$ be a finite connected $d$-regular graph with~$d \geq 3$, endowed with nonnegative edge lengths such that the subgraph formed by the zero-length edges contains no cycle.
  Let $x$ be a vertex of~$X$.
  Denote by~$\tilde{x}$ a lift of~$x$ to the universal cover~$\tilde{X}$.
  Then, for every $r>0$,
  \[
    \length \, B(\tilde{x},r) = dr + (d-2) \, \sum_{k \geq 1} \sum_{\gamma \in \mathcal{P}_k(x)} \left( r - \ell(\gamma) \right)_+ .
  \]
\end{lemma}

\begin{proof}
  The universal cover $\tilde{X}$ is a $d$-regular  tree equipped with the pseudometric induced by the edge lengths.
  Since the subgraph formed by the zero-length edges contains no cycle, its connected components lift to~$\tilde X$ and only finitely many vertices lie at distance at most~$r$ from~$\tilde x$.

  For almost every $t>0$, the distance sphere~$S(\tilde x,t)$ has
  \[
    |S(\tilde{x},t)| = d + (d-2) \, | \{ \tilde{y} \textrm{ vertex of } \tilde{X} \textrm{ with } \tilde{y} \neq \tilde{x} \textrm{ and }  \dist(\tilde{x},\tilde{y}) < t \} |
  \]
  points.
  Indeed, whenever the sphere passes through a vertex, one incoming branch is replaced by $d-1$ outgoing branches,
  increasing its cardinality by $d-2$.
  This remains valid when the sphere passes through vertices at zero distance from each other.

  By the coarea formula, integrating with respect to~$t$ gives
  \begin{align*}
    \length \, B(\tilde{x},r) & = \int_0^r |S(\tilde{x},t)| \, dt                                                                    \\
                              & = dr + (d-2) \, \sum_{\tilde{y} \neq \tilde{x}} \left( r - \dist(\tilde{x},\tilde{y}) \right)_+      \\
                              & = dr + (d-2) \, \sum_{k \geq 1} \sum_{\gamma \in \mathcal{P}_k(x)} \left( r - \ell(\gamma) \right)_+
  \end{align*}
  where the last equality follows from the correspondence between the vertices~\mbox{$\tilde{y} \neq \tilde{x}$} and the edge-paths $\gamma \in \bigsqcup_{k \geq 1} \mathcal{P}_k(x)$ joining~$x$ to~$y$, with $\dist(\tilde{x},\tilde{y})=\ell(\gamma)$; see Remark~\ref{rem:correspondence}.
\end{proof}

\begin{example}
  Let $T_d$ be the $d$-regular tree with all the edges of unit length.
  By Lemma~\ref{lem:ball_paths} and the relation~\eqref{eq:Pkx}, the length of a ball of radius~$r$ centered at a vertex is
  \begin{equation} \label{eq:Vd}
    V_{T_d}(r) = dr + d(d-2)  \, \sum_{k \geq 1} (d-1)^{k-1} ( r - k )_+.
  \end{equation}
\end{example}

We conclude this section with an alternative expression for the length of a ball in the universal cover, obtained by summing its intersections with individual edges.

\begin{lemma}\label{lem:ball_all}
  Let $X$ be a finite connected graph endowed with a length metric.
  Let $x$ be a vertex of~$X$.
  Then, for every $r>0$,
  \begin{equation}\label{eq:edgesum}
    \length\,B(\tilde{x},r)=\sum_{\gamma a} \min\{\ell(a),(r-\ell(\gamma))_+\},
  \end{equation}
  where the sum runs over all nonempty oriented edge-paths~$\gamma a$ starting at~$x$.
  Only finitely many terms are nonzero.
\end{lemma}

\begin{proof}
  By Remark~\ref{rem:correspondence}, the vertices of~$\tilde{X}$ can be identified with the oriented edge-paths starting at~$x$, with $\emptyset_x$ corresponding to $\tilde{x}$.
  Each edge of~$\tilde{X}$ then joins~$\gamma$ to~$\gamma a$ and has length~$\ell(a)$.
  Since its points lie at distances between~$\ell(\gamma)$ and~$\ell(\gamma)+\ell(a)$ from~$\tilde x$, its intersection with~$B(\tilde{x},r)$~has length
  \[
    \min\{\ell(a),(r-\ell(\gamma))_+\}.
  \]
  Summing over all edges gives~\eqref{eq:edgesum}.
  Finally, since every edge of~$X$ has positive length, only finitely many paths~$\gamma$ satisfy~$\ell(\gamma)<r$.
\end{proof}

We also record the effect of rescaling the edge lengths.

\begin{remark} \label{rem:scaling}
  Let $X$ be a finite connected graph endowed with a length metric.
  For $\lambda>0$, denote by~$X_\lambda$ the graph~$X$ with all its edge lengths multiplied by~$\lambda$.
  In particular,
  \[
    \length(X_\lambda)=\lambda \, \length(X).
  \]
  Moreover, the universal cover~$\tilde{X}_\lambda$ of~$X_\lambda$ is~$\tilde{X}$ with all its edge lengths multiplied by~$\lambda$.
  Thus, for every vertex~$x$ and every $r>0$,
  \[
    \length \, B_{\tilde{X}_\lambda}(\tilde{x},\lambda r)= \lambda \, \length \, B_{\tilde{X}}(\tilde{x},r).
  \]
\end{remark}

\section{Regular graphs: proof of Theorem~\ref{theo:main}}\label{sec:proof13}

In this section, we prove Theorem~\ref{theo:main} by averaging the lengths of edge-paths and applying a convexity argument.
The proof also applies for nonnegative edge lengths, provided that the subgraph formed by the zero-length edges contains no cycle.
We also characterize the radii at which equality holds.

\begin{proof}[Proof of Theorem~\ref{theo:main}]
  Let $X$ be a finite connected $d$-regular graph with~$d \geq 3$, endowed with nonnegative edge lengths such that the subgraph formed by the zero-length edges contains no cycle and $\length(X) \leq E(X)$.

  Averaging the lengths of the balls~$B(\tilde{x},r)$ over the vertices of~$X$ and using Lemma~\ref{lem:ball_paths}, we obtain
  \begin{equation} \label{eq:average}
    \frac{1}{V(X)} \sum_{x \in V(X)} \length \, B(\tilde{x},r) = dr + \frac{d-2}{V(X)} \sum_{k \geq 1} \sum_{\gamma \in \mathcal{P}_k} \left( r - \ell(\gamma) \right)_+.
  \end{equation}
  The sum of the right-hand side is finite for every fixed~$r$, since the subgraph formed by the zero-length edges contains no cycle.

  Set $\lambda = \frac{\length(X)}{E(X)}$.
  Applying Jensen's inequality to the convex function~$\varphi_r(s) = (r-s)_+$ yields
  \begin{equation} \label{eq:jensen}
    \frac{1}{|\mathcal{P}_k|} \sum_{\gamma \in \mathcal{P}_k} \left( r - \ell(\gamma) \right)_+ \geq \left( r - \frac{1}{|\mathcal{P}_k|} \sum_{\gamma \in \mathcal{P}_k} \ell(\gamma) \right)_+ = (r- \lambda k)_+ \geq (r-k)_+
  \end{equation}
  where the equality follows from Lemma~\ref{lem:w-comb} and the last inequality from $\length(X) \leq E(X)$.

  Substituting into~\eqref{eq:average}, we obtain
  \begin{align}
    \frac{1}{V(X)} \sum_{x \in V(X)} \length \, B(\tilde{x},r) & \geq dr + (d-2) \, \sum_{k \geq 1}  \frac{|\mathcal{P}_k|}{V(X)} \, (r-k)_+ \nonumber \\
                                                               & = dr + d(d-2) \, \sum_{k \geq 1} (d-1)^{k-1} (r-k)_+ \nonumber                        \\
                                                               & = V_{T_d}(r) \label{eq:average2}
  \end{align}
  where the first equality follows from~\eqref{eq:Pk}, and the second from~\eqref{eq:Vd}.

  The lower bound~\eqref{eq:main} on the maximal length of an $r$-ball of~$\tilde{X}$ immediately follows.

  Suppose now that equality holds in the averaged inequality~\eqref{eq:main-average} or in the maximal ball length bound~\eqref{eq:main} at some integer radius~$r=n \geq 1$.
  Since the average ball length is at least~$V_{T_d}(n)$, equality also holds in the averaged inequality~\eqref{eq:average2}.
  In particular, equality holds in Jensen's inequality~\eqref{eq:jensen} for~$k=n$.

  For $r=n=k$, the right-hand side of~\eqref{eq:jensen} is equal to zero.
  Therefore, equality in~\eqref{eq:jensen} implies that $\ell(\gamma) \geq n$ for every $\gamma \in \mathcal{P}_n$.
  Furthermore, it also implies that $\lambda=1$.
  That is, $\length(X)=E(X)$.
  On the other hand, Lemma~\ref{lem:w-comb} shows that the average of the lengths~$\ell(\gamma)$ with $\gamma \in \mathcal{P}_n$ is equal to~$n$.
  Hence, $\ell(\gamma)=n$ for every $\gamma \in \mathcal{P}_n$.

  We claim that all edges have the same length.
  If $n=1$, this follows immediately from the preceding equality.
  Suppose that $n \geq 2$.
  Fix a vertex~$x$ and two distinct oriented edges~$a$ and~$b$ issuing from~$x$.
  Since $d\ge3$, there exists a third oriented edge~$c$ issuing from~$x$, distinct from both~$a$ and~$b$.
  Choose an oriented edge-path~$\gamma$ of combinatorial length~$n-1$ ending at~$x$ with last edge~$\bar{c}$. Then $\gamma a$ and~$\gamma b$ both belong to~$\mathcal P_n$, so
  \[
    \ell(\gamma) + \ell(a) = \ell(\gamma a) = n = \ell(\gamma b) = \ell(\gamma) + \ell(b).
  \]
  Thus, $\ell(a)=\ell(b)$.
  Since $a$ and~$b$ were arbitrary, all edges incident to~$x$ have the same length.
  By connectedness, all edges of~$X$ have the same length.
  The normalization $\length(X)=E(X)$ then implies that every edge has length~$1$.

  Conversely, if all edges have length~$1$, the universal cover is isometric to~$T_d$, and equality holds for every radius~$r>0$.
\end{proof}

The following result gives a characterization of the equality case in Theorem~\ref{theo:main}.

\begin{proposition} \label{prop:equality}
  Let $X$ be a finite connected $d$-regular graph with~$d \geq 3$, endowed with a length metric normalized so that $\length(X) = E(X)$.
  For every integer $k\ge1$, denote by~$m_k$ and~$M_k$ the smallest and the largest lengths of the oriented edge-paths in~$\mathcal{P}_k$, and set $M_0=0$.
  Then, for every~\mbox{$r>0$},
  \begin{equation} \label{eq:equality}
    \frac{1}{V(X)} \sum_{x \in V(X)} \length \, B(\tilde{x},r) = V_{T_d}(r)
  \end{equation}
  if and only if $M_k \leq r \leq m_{k+1}$ for some integer $k \geq 0$.
\end{proposition}

\begin{proof}
  Recall that the average length of the paths in~$\mathcal P_k$ is equal to~$k$, by Lemma~\ref{lem:w-comb}.
  Since the function $\varphi_r(s)=(r-s)_+$ is affine on $(-\infty,r]$ and on $[r,\infty)$, equality holds in Jensen's inequality~\eqref{eq:jensen} if and only if all the paths in~$\mathcal{P}_k$ have length at most~$r$, or all have length at least~$r$.
  This condition can be written as
  \begin{equation} \label{eq:condition}
    M_k \leq r \qquad \textrm{ or } \qquad r \leq m_k.
  \end{equation}

  By the expressions~\eqref{eq:average} and~\eqref{eq:Vd}, along with the formula~\eqref{eq:Pk}, equality holds in~\eqref{eq:equality} if and only if equality holds in~\eqref{eq:jensen} for every~$k \geq 1$, which means exactly that the condition~\eqref{eq:condition} is satisfied for every~$k \geq 1$.

  Since all edge lengths are positive, both sequences~$(m_k)$ and~$(M_k)$ are increasing and tend to infinity.

  Fix $r>0$.
  Since $M_0=0$ and $(M_k)$ tends to infinity, there is a largest integer $K \geq 0$ with $M_K \leq r$.
  For every $k \leq K$, we have $M_k \leq r$, so the condition~\eqref{eq:condition} is satisfied for all such~$k$.
  For every $k \geq K+1$, we have $M_k > r$.
  Hence, the condition~\eqref{eq:condition} is satisfied for all such~$k$ if and only if $r \leq m_k$ for every $k \geq K+1$.
  Since $(m_k)$ is increasing, this is equivalent to~$r \leq m_{K+1}$.

  Therefore, equality holds in~\eqref{eq:equality} if and only if
  \[
    M_K \leq r \leq m_{K+1}.
  \]
  This proves one implication of Proposition~\ref{prop:equality}.

  If $M_k \leq r \leq m_{k+1}$ for some integer~$k \geq 0$, then $k \leq K$.
  Suppose that $k <K$.
  Then $r \leq m_{k+1} \leq m_K \leq M_K \leq r$.
  Thus, $m_K=r<m_{K+1}$.
  Hence, $M_K \leq r \leq m_{K+1}$, which implies~\eqref{eq:equality}.
\end{proof}

We apply the previous result in the following example to show that rigidity fails at every noninteger radius.

\begin{example} \label{ex:noninteger}
  Fix a noninteger radius~$r>0$.
  Let $k=\lfloor r \rfloor$.
  Note that $k < r < k+1$.
  Choose $\delta >0$ with $k(1+\delta) <r$ and $r < (k+1)(1-\delta)$.
  Let $X$ be a finite connected $d$-regular graph with~$d \geq 3$.
  Endow~$X$ with a nonequilateral metric, where two edges have lengths $1+\delta$ and~$1-\delta$, and all the other edges have unit length.
  By construction, $\length(X)=E(X)$ and $M_k < r < m_{k+1}$.
  Thus, equality holds in~\eqref{eq:main-average} by Proposition~\ref{prop:equality}.

  If $X$ is the $d$-dipole formed by two vertices joined by $d$ edges, then equality also holds in~\eqref{eq:main}.
  Indeed, the isometry of~$X$ exchanging its two vertices shows that the $r$-balls centered at the two vertices have the same length, so both are equal to their average~$V_{T_d}(r)$.
\end{example}

\section{Arbitrary graphs: proof of Theorem~\ref{theo:nonreg}} \label{sec:nonreg}

In this section, we prove Theorem~\ref{theo:nonreg} by reducing the general case to the trivalent case of Theorem~\ref{theo:main}, replacing vertices of higher degree by trivalent trees with zero-length edges.

\begin{proof}[Proof of Theorem~\ref{theo:nonreg}]
  Let $X$ be a finite connected graph with cyclomatic number~$c=c(X) \geq 2$, endowed with a length metric satisfying
  \[
    \length(X) \leq 3(c-1).
  \]
  By construction of the core~$K$ of~$X$, we have
  \[
    c(K) =c(X) \qquad \textrm{ and } \qquad \length(K) \leq \length(X).
  \]

  For each vertex $x \in K$ of degree~$q \geq 4$, replace~$x$ by a trivalent tree with $q$ leaves, attach the $q$ original incident edges to these leaves, assign length zero to all the new edges, and suppress the resulting degree-two vertices, so that every edge of~$K$ keeps its length.
  The resulting graph~$Z$ is trivalent.

  Since a trivalent tree with $q$ leaves has $q-2$ trivalent vertices and $q-3$ internal edges, we obtain
  \[
    V(Z) = \sum_{x \in V(K)} \left( \deg_K(x)-2 \right) = 2 E(K) - 2 V(K) = 2(c-1).
  \]
  Consequently,
  \[
    E(Z) = 3(c-1).
  \]
  Moreover,
  \[
    \length(Z) = \length(K) \leq 3(c-1)=E(Z).
  \]
  The subgraph of~$Z$ formed by the edges of zero length contains no cycle, so the comparison theorem for regular graphs applies to~$Z$.

  By Theorem~\ref{theo:main},
  \[
    \frac{1}{2(c-1)} \sum_{z \in V(Z)} \length \, B_{\tilde{Z}}(\tilde{z},r) \geq V_{T_3}(r).
  \]

  Each vertex $x \in V(K)$ corresponds to exactly $\deg_K(x)-2$ vertices of~$Z$.
  Collapsing each lift of each new tree to a point maps~$\tilde Z$ onto~$\tilde K$, sending the lifts of the vertices~$z$ corresponding to~$x$ to lifts of~$x$.
  Since the collapsed edges have length zero and the other edges keep their lengths, this map preserves distances, and is one-to-one and length-preserving outside the collapsed edges.
  Therefore,
  \[
    \length \, B_{\tilde{Z}}(\tilde{z},r) = \length \, B_{\tilde{K}}(\tilde{x},r)
  \]
  for every vertex $z$ corresponding to $x$.
  It follows that
  \begin{equation*} \label{eq:sumK}
    \frac{1}{2(c-1)} \sum_{x \in V(K)} \left( \deg_K(x)-2 \right)  \length \, B_{\tilde{K}}(\tilde{x},r) \geq V_{T_3}(r)
  \end{equation*}
  which proves~\eqref{eq:nonreg-average}.

  In particular, there exists a vertex~$x \in K$ such that
  \[
    \length \, B_{\tilde{X}}(\tilde{x},r) \geq \length \, B_{\tilde{K}}(\tilde{x},r) \geq V_{T_3}(r).
  \]
  This proves~\eqref{eq:nonreg}.

  Suppose now that equality holds in either~\eqref{eq:nonreg-average} or~\eqref{eq:nonreg} at some integer radius~$r=n \geq 1$.
  Since the average ball length is at least~$V_{T_3}(n)$, equality also holds in the averaged inequality~\eqref{eq:nonreg-average}.
  By construction, the left-hand side of~\eqref{eq:nonreg-average} is equal to the average length of the balls of radius~$n$ centered at lifts of vertices of~$Z$.
  Since the subgraph of~$Z$ formed by the edges of zero length contains no cycle, Theorem~\ref{theo:main}, which also applies to nonnegative edge lengths, implies that every edge of~$Z$ has unit length.
  In particular, $Z$ has no zero-length edges.
  Hence, $K=Z$ is trivalent with all edges of unit length, and $\length(K)=E(K)=3(c-1)$.
  Since $\length(K) \leq \length(X) \leq 3(c-1)$, we have $\length(K)=\length(X)$.
  Thus, no edges were removed in the construction of the core, and $X$ is isometric to the equilateral trivalent graph~$K$.

  Conversely, if $X$ is isometric to an equilateral trivalent graph with unit edge lengths, then its universal cover is isometric to~$T_3$, and equality holds in both~\eqref{eq:nonreg-average} and~\eqref{eq:nonreg} at every radius.
\end{proof}

\section{Non-rigidity: theta graphs}\label{sec:theta}

We show that, unlike in the surface case~\cite{CV}, equality in the optimal ball length comparison~\eqref{eq:main} does not imply rigidity, even when the center is allowed to vary over the entire universal cover.

\medskip

The \emph{theta graph}~$\Theta$ consists of two vertices joined by three edges.
It is $3$-regular, with $V(\Theta)=2$, $E(\Theta)=3$ and $c(\Theta)=2$.

\begin{proposition}\label{prop:theta}
  Let $\Theta$ be the theta graph endowed with the length metric whose edges have lengths $x < y < z <\frac{3}{2}$, with $x+y+z=3$, so that $\length(\Theta)=E(\Theta)$.
  Then, for every $r\in[z,x+y]$,
  \[
    \max_{p \in \tilde{\Theta}} \length \, B(p,r) = V_{T_3}(r)
  \]
  where the maximum is over all the points of~$\tilde\Theta$.

  Moreover, the maximum is attained precisely at the vertices of~$\tilde\Theta$.
  In particular, equality holds in Theorem~\ref{theo:main} for every $r\in[z,x+y]$.
\end{proposition}

\begin{proof}
  Since $x < y <  z <\frac{3}{2}$ and $x+y+z=3$, we have $1 < z < x+y < 2$.
  Thus, $[z,x+y] \subseteq (1,2)$.
  Fix $r\in[z,x+y]$.

  Let us compute the length of~$B(\tilde{v},r)$, where $v$ is a vertex of~$\Theta$.
  The three oriented edge-paths of~$\mathcal{P}_1(v)$ have lengths $x,y,z\le r$, while every oriented edge-path with at least two edges has length at least $x+y\ge r$.
  By Lemma~\ref{lem:ball_paths} with $d=3$, it follows that
  \[
    \length \, B(\tilde v,r)=3r+(r-x)+(r-y)+(r-z)=6r-3 = V_{T_3}(r).
  \]

  Let us show that this is the maximal length of a ball~$B(p,r)$, where $p \in \tilde{\Theta}$.
  The point~$p$ lies on an edge~$e=(v,v')$ of length~$\ell$ with endpoints~$v$ and~$v'$.
  Write $s=\dist(p,v)$ and $s'=\dist(p,v')$, so that $s+s'=\ell$.
  Since $\ell \leq z \leq r$, the whole edge~$e$ lies in~$B(p,r)$.

  Let $\ell_1$ and~$\ell_2$ be the other two edge lengths of~$\Theta$.
  The edge~$e$ is adjacent to two edges~$(v,v_1)$ and~$(v,v_2)$ of length~$\ell_1$ and~$\ell_2$ at~$v$.
  Since $r \leq x+y$, the ball can meet the two edges beyond~$v_i$ along segments, each of length exactly~$(r-s-\ell_i)_+$, but cannot go beyond those two edges.
  Thus, we have
  \[
    \length \left( B(p,r) \cap T_i \right)= \min \{ r-s, \ell_i \} + 2 (r-s-\ell_i)_+ = (r-s) + (r-s-\ell_i)_+
  \]
  where $T_i$ is the tripod consisting of the edge~$(v,v_i)$ and the two edges issuing from~$v_i$ other than~$(v,v_i)$.
  The same formula holds at~$v'$ with~$s'$ in place of~$s$.
  Summing the different parts gives
  \[
    \length \, B(p,r) = 4r - \ell + \sum_{i=1}^2 \left[ (r-s-\ell_i)_+ + (r-s'-\ell_i)_+ \right].
  \]

  Let $r_i=r-\ell_i$.
  Since $\ell_i \leq z \leq r \leq x+y \leq \ell+\ell_i$, we have $0 \leq r_i \leq \ell$.
  If both $r_i-s$ and $r_i-s'$ are positive, their sum is $2r_i-\ell \leq r_i$; if only one is, it is at most~$r_i$; and if neither is positive, the inequality is immediate.
  Thus,
  \[
    (r_i-s)_+ + (r_i-s')_+ \leq r_i
  \]
  with equality if and only if $s=0$, $s'=0$, $r_i=0$ or $r_i=\ell$.
  It follows that
  \begin{equation} \label{eq:theta}
    \length \, B(p,r) \leq 4r-\ell+(r-\ell_1)+(r-\ell_2)=6r-3 = V_{T_3}(r)
  \end{equation}
  where the first equality comes from $\ell+\ell_1+\ell_2=x+y+z=3$.

  It remains to determine when equality holds.
  Suppose that $p$ lies in the interior of~$e$, so that $s,s'>0$.
  By the previous analysis, equality in~\eqref{eq:theta} holds if and only if $r \in \{ \ell_1, \ell_1+\ell \} \cap \{ \ell_2, \ell_2+\ell \}$.
  Since the edges of~$\Theta$ have different lengths and each edge length is smaller than the sum of the other two (because $z<x+y$), this intersection is empty.
  Thus, equality in~\eqref{eq:theta} holds precisely when $p$ is a vertex of~$\tilde{\Theta}$.
\end{proof}

\begin{remark}
  The equality for balls centered at vertices in Proposition~\ref{prop:theta} also follows from Proposition~\ref{prop:equality}.
  Indeed, for the theta graph, $M_1=z$ and $m_2=x+y$, since consecutive edges of an edge-path are distinct.
  Thus, $[z,x+y]=[M_1,m_2]$.
  Since the map exchanging the two vertices of~$\Theta$ is an isometry, both vertex balls are equal to their average.
  Proposition~\ref{prop:theta} strengthens this equality by showing that no ball centered at an arbitrary point of~$\tilde{\Theta}$ has
  greater length.
\end{remark}

\end{document}